\documentclass[12pt]{article}
\title{Real embeddings and the Atiyah-Patodi-Singer index theorem for
Dirac operators\thanks{Partially 
supported by NSF grant DMS-9022140 
when both authors were visiting MSRI in 1994.}}
\author{Xianzhe Dai\thanks{Partially supported by NSF and 
Alfred P. Sloan Foundation.} $\ $ and Weiping Zhang\thanks{Partially 
supported by  the CNNSF,
EMC  and the Qiu Shi Foundation.}}
\begin{document}
\maketitle
\begin{abstract} We present  the details of our embedding proof, 
 which was announced in [DZ1],
 of the Atiyah-Patodi-Singer index theorem for Dirac operators on manifolds
with boundary [APS1]. 
\end{abstract}

{\bf Introduction}

$\ $

The index theorem of Atiyah, Patodi and Singer [APS1, (4.3)] for Dirac operators on
manifolds with boundary has played important roles in various problems in geometry,
topology as well as  mathematical physics. Not surprisingly then, there are by now
quite a number of proofs of this index theorem other than Atiyah, Patodi and Singer's
original proof  [APS1]. Among these proofs
we mention those of Cheeger [C1, 2] (see also Chou [Ch]), Bismut-Cheeger [BC1] and Melrose [M].
One common point underlying all these proofs (including the original one) is that they can all be viewed, in one way or another,
as certain extensions to 
manifolds with boundary of the heat kernel proof of the local index theorem 
for Dirac operators on closed manifolds (cf. [BeGV]). That is, one starts with
a Mckean-Singer type formula and then studies the small time asymptotics of the 
corresponding heat kernels. In particular, one makes use of
the explicit formulas for the heat kernel of the Laplace operators
on the cylinder ([APS1], [M]) and/or cone ([BC1], [C1, 2], [Ch]) 
(being attached the boundary)
for the analysis near the boundary. The $\eta$-invariant on the boundary, which
was first defined in [APS1], appears naturally during the
process.

Now recall that Atiyah and Singer [AS] also have a $K$-theoretic proof of 
their index theorem for elliptic operators on closed manifolds. In such a proof,
one transforms the problem, through direct image constructions in $K$-theory, to
a sphere and then applies the Bott periodicity theorem on the sphere 
to establish the result.\footnote{See also the book of Lawson-Michelsohn [LM] for a
comprehensive treatment of this approach.}
It is thus natural to ask whether the strategy of Atiyah-Singer's $K$-theoretic ideas
can be used to prove the Atiyah-Patodi-Singer index theorem for manifolds with
boundary. The purpose of this paper is to present such a proof, of which an
announcement of basic ideas has already appeared in [DZ1].

Briefly speaking, we  embed the manifold with boundary 
under consideration into a ball, instead of
a sphere, so that it maps the boundary of the original manifold to the boundary
sphere of the ball, and reduce the problem to the ball. Now since any
vector bundle on the ball is topologically trivial, one obtains the result immediately.
 This works even when the original manifold has no boundary, 
giving a proof of the Atiyah-Singer index theorem for Dirac operators. The Bott
periodicity theorem is thus not needed.

Observe that in [AS], Atiyah and Singer made heavy use of the techniques of
pseudodifferential operators, which is not suitable for treating directly the 
global elliptic boundary problems. This is the first serious difficulty in extending directly
the arguments in [AS] to deal with the Atiyah-Patodi-Singer boundary problems. 

On the other hand, Bismut and Lebeau developed in [BL] a general and direct
localization procedure which applies to a wide range of localization problems involving Dirac type operators. 
For example, it has lead to a direct analytic
treatment of the index theorem for Dirac operators on closed manifolds along the
lines of [AS] (cf. [Z, Remark 2.6]), as well as a localization
formula for $\eta$-invariants of Dirac operators [BZ] which may be viewed as an odd 
dimensional analogue of the main result in [BL]. It is these techniques and results
that  will be used in the present paper, giving an embedding  proof of the
Atiyah-Patodi-Singer index theorem for {\it Dirac} operators on manifolds with
boundary [APS1].

In the proof described in [DZ1], we also used in an essential way
Cheeger's cone method [C1, 2]. The reason being, in order to apply
Bismut-Lebeau's method [BL], 
we need to transfer the Atiyah-Patodi-Singer boundary problem to an elliptic problem
on  certain manifolds with cone-like singularity. Now, in the present paper, we will show that how one can avoid the analysis on  the cone at all. This is done by considering the 
Atiyah-Patodi-Singer type boundary value problem for certain non-differetial
operators arising naturally from the analysis in [BL]. In this way, one no longer
encounters the heat kernel analysis on cylinders and/or cones which are
essential for the other proofs of the Atiyah-Patodi-Singer index theorem.
We regard this as a major technical simplification with respect to [DZ1].

In a separate paper [DZ3], we will further 
extend the main result of this paper to  the case of families.
In particular, we will give a new proof of the family index theorem 
of Bismut-Cheeger [BC1, 2] and Melrose-Piazza [MP] along the lines of this paper.

This paper is organized as follows. In Section 1, we prove an important 
variation formula for the indices of the Atiyah-Patodi-Singer boundary value
problems for Dirac operators on manifolds with boundary. In Section 2, 
we state a localization formula of Riemann-Roch type for the indices of 
the Atiyah-Patodi-Singer boundary value
problems for Dirac operators on manifolds with boundary. 
In Section 3, we prove the Riemann-Roch property stated in Section 2.
In Section 4, by combining the results in Sections 1, 2 with those of
Bismut-Zhang [BZ], we complete our proof of the Atiyah-Patodi-Singer index theorem
for Dirac operators on manifolds with boundary. There is also an appendix in which 
we prove a harmonic oscillator property for certain
Dirac operators on flat spaces, which
plays an essential role in the main text.

$$\ $$

{\bf \S 1. Dirac type operators on manifolds with boundary: index and its
variations}

$\ $

In this section, we recall the definition of the Atiyah-Patodi-Singer boundary 
value problems
[APS1] for Dirac type operators on {\it Spin} manifolds with boundary.
We also prove an important variation formula for the indices of these boundary value problems.

$\ $

Let $X$ be a compact
 oriented even dimensional spin manifold with boundary $\partial X$.
We assume that $X$ has been equipped with a fixed spin structure. Then $\partial X$
carries the canonically  induced orientation and spin structure.

Let $g^{TX}$ be a metric on $TX$. Let $g^{T\partial X}$ be its restriction on
$T\partial X$. We assume that $g^{TX}$ is of product structure near
the boundary $\partial X$. That is, there is an open neighborhood 
$U_\alpha =[0,\alpha)\times \partial X$ of $\partial X$ in $X$ with $\alpha>0$
such that one has the orthogonal splitting on $U_\alpha$,
$$\left.g^{TX}\right|_{U_\alpha}=dr^2\oplus \pi^*_\alpha g^{T\partial X},\eqno (1.1)$$
where $\pi_\alpha
:[0,\alpha)\times \partial X\rightarrow \partial X$ is the obvious projection
onto the second factor.

Let $\xi$ be a Hermitian vector bundle over $X$ with Hermitian metric $g^\xi$.
Let $\nabla^{\xi}$ be a Hermitian connection on $\xi$ with respect to $g^\xi$.
We make the assumption that over the open neighborhood $U_\alpha$  of $\partial X$,
one has 
$$\left. g^\xi\right|_{U_\alpha}=\pi^*_\alpha\left(g^\xi|_{\partial X}\right),\ \ 
\left.\nabla^\xi\right|_{U_\alpha}=\pi^*_\alpha\left(\nabla^\xi |_{\partial X}\right).
\eqno (1.2)$$

By taking $\alpha>0$ sufficiently small, one can always find $g^{TX}$, $g^\xi$ and
$\nabla^\xi$ verifying (1.1) and (1.2).

Let $S(TX)=S_+(TX)\oplus S_-(TX)$ be the ${\bf Z}_2$-graded Hermitian vector
bundle of spinors associated to $(TX, g^{TX})$. Let $\nabla^{S(TX)}$
be the Hermitian connection on $S(TX)$ canonically induced from the Levi-Civita
connection $\nabla^{TX}$ of $g^{TX}$. Then $\nabla^{S(TX)}$ preserves the 
${\bf Z}_2$-splitting $S(TX)=S_+(TX)\oplus S_-(TX)$. We denote by $\nabla^{S_\pm (TX)}$
the restriction of $\nabla^{S(TX)}$ on $S_\pm (TX)$. Let 
$\nabla^{S(TX)\otimes \xi}$ (resp. $\nabla^{S_\pm (TX)\otimes \xi}$) be the Hermitian 
connection on $S(TX)\otimes \xi$ (resp. $S_\pm (TX)\otimes \xi$) 
obtained from the tensor product of $\nabla^{S(TX)}$ (resp. $\nabla^{S_\pm (TX)}$)
and $\nabla^\xi$.

For any $e\in TX$, let $c(e)$ be the Clifford action of $e$ on $S(TX)$.
Then $c(e)$ extends to an action on $S(TX)\otimes \xi$ by acting as identity on $\xi$.
We still denote this extended action by $c(e)$.

Let $e_1,\cdots,e_{\dim X}$ be an oriented (local) orthonormal base of $TX$. We can then
define the (total) twisted Dirac operator with 
coefficient bundle $\xi$ as follows (cf. [BeGV] and [LM]),
$$D^\xi =\sum_{i=1}^{\dim X}c(e_i)\nabla^{S(TX)\otimes \xi}_{e_i}:
\Gamma (S(TX)\otimes \xi) \rightarrow \Gamma (S(TX)\otimes \xi) .\eqno (1.3)$$
Let $D^\xi_\pm$ be the restriction of $D^\xi$ on $\Gamma(S_\pm (TX)\otimes \xi)$. Then
$D^\xi_-$ is the formal adjoint of $D^\xi_+$.

$\ $

{\bf Definition 1.1.}  {\it By a Dirac type operator on $\Gamma (S(TX)\otimes \xi) $, 
 we mean a first order differential operator 
$D:\Gamma (S(TX)\otimes \xi) \rightarrow \Gamma (S(TX)\otimes \xi) $ such that
$D-D^\xi$ is an odd self-adjoint element
of zeroth order, and that for $\alpha>0$ sufficiently small, the following 
identity holds on $U_\alpha$,
$$D=c\left({\partial \over \partial r}\right)\left({\partial \over \partial r}+B
\right), \eqno (1.4)$$
with $B$ independent of $r$ and its restriction on
$\Gamma (S(TX)\otimes \xi)|_{\partial X}$ formally self-adjoint. We will also
call the  restriction $D_+$ (resp. $D_-$) of $D$ to
$\Gamma (S_+(TX)\otimes \xi) $ (resp. $\Gamma (S_-(TX)\otimes \xi)$)
a Dirac type operator.}

$\ $

When there is no confusion, we will also use $B$ to denote its restriction on 
$(S(TX)\otimes \xi)|_{\partial X}$. Clearly, $B$ preserves the 
${\bf Z}_2$-grading of $(S(TX)\otimes \xi)|_{\partial X} =
 (S_+(TX)\otimes \xi)|_{\partial X} \oplus (S_-(TX)\otimes \xi)|_{\partial X} $. 
We denote by $B_\pm$ the restriction of $B$ on $(S_\pm (TX)\otimes \xi)|_{\partial X}$.

Now consider the formally self-adjoint first order differential
operator $B_+$, which is clearly elliptic,  acting on sections of
$(S_+ (TX)\otimes \xi)|_{\partial X}$. Then the $L^2$-completion of 
$(S_+ (TX)\otimes \xi)|_{\partial X} $ admits an orthogonal decomposition
$$L^2\left((S_+ (TX)\otimes \xi)|_{\partial X}\right) =
\bigoplus_{\lambda \in  {\rm Spec}(B_+)} E_\lambda  ,\eqno (1.5)$$
where $E_\lambda$ is the eigenspace of $\lambda$. 

For any $a\in {\bf R}$, let  $L^2_{\geq a} ((S_+ (TX)\otimes \xi)|_{\partial X}) $ 
denote the direct sum of the eigenspaces $E_\lambda$ associated to the 
eigenvalues $\lambda \geq a$. 
Let $P_{+,\geq a} $ denote the orthogonal projection from 
$L^2((S_+ (TX)\otimes \xi)|_{\partial X}) $ to 
$L^2_{\geq a}((S_+ (TX)\otimes \xi)|_{\partial X}) $.
We  call the particular projection $P_{+,\geq 0} $ the Atiyah-Patodi-Singer projection
associated to $B_+$, to emphasize its role in [APS1].

Following [APS1], one can then impose the boundary value problem
$$\left(D_+, P_{+,\geq a} \right) :\left\{ u:u\in \Gamma (S_+(TX)\otimes \xi) ,\ 
P_{+,\geq a} u|_{\partial X}=0\right\}\rightarrow \Gamma (S_-(TX)\otimes \xi),
\eqno (1.6)$$
which is Fredholm by [APS1]. In particular, we call the boundary problem
$(D_+,P_{+,\geq 0})$ the Atiyah-Patodi-Singer boundary problem associated
to $D_+$. We denote by ${\rm ind}(D_+, P_{+,\geq a}) $ the index of 
the associated Fredholm operator.

Now let $D_+(s)$, $0\leq s\leq 1$, be a smooth family of Dirac type operators
with the induced boundary operators $B_+(s)$.
We can now state the main result of this section,
which has been announced in [DZ1, Theorem 1.1],  as follows.

$\ $

{\bf Theorem 1.2.} {\it The following identity holds,
$${\rm ind}\left(D_+(1), P_{+,\geq 0}(1)\right) - 
{\rm ind}\left(D_+(0), P_{+,\geq 0}(0)\right) 
=-{\rm sf}\{ B_+(s), 0\leq s\leq 1\},\eqno (1.7)$$
where ${\rm sf}$ is the notation for the spectral flow of Atiyah-Patodi-Singer
[APS2].}

{\it Proof.}  Take any $0\leq s_0\leq 1$. Let $2a_0$ be the 
minimal absolute value of the
nonzero eigenvalues of $B_+(s_0)$. Then there exsists $\varepsilon_0>0$ such that
for any $s\in [s_0-\varepsilon_0 ,s_0+\varepsilon_0 ]\cap [0,1]$, $a_0$ is
not an eigenvalue of $B_+(s)$. Then for any $s\in [s_0-\varepsilon , s_0+\varepsilon_0 ]\cap [0,1]$, $(D_+(s), P_{+,\geq a_0}(s))$ defines a continuous family of Fredholm operators. Therefore,
$$ {\rm ind}(D_+(s), P_{+,\geq a_0}(s)) = {\rm ind}(D_+(s_0), P_{+,\geq a_0}(s_0)) .
\eqno (1.8)$$

On the other hand, by the classical Agranovi\v{c}-Dynin type formula
(cf. [BoW, Chap. 21]) and the definition of spectral flow [APS2], 
one verifies easily that
$$ {\rm ind}(D_+(s), P_{+,\geq a_0}(s)) -{\rm ind}(D_+(s), P_{+,\geq 0}(s)) 
=-{\rm sf}\{B_+(s)+ua_0,0\leq u\leq 1\},$$
$$ {\rm ind}(D_+(s_0), P_{+,\geq a_0}(s_0)) -{\rm ind}(D_+(s_0), P_{+,\geq 0}(s_0))$$ 
$$=-{\rm sf}\{B_+(s_0)+ua_0,0\leq u\leq 1\}.\eqno (1.9)$$

Formula (1.7) follows easily from (1.8), (1.9) and 
the additivity (using twice) of the spectral flow [APS2]. Q.E.D.

$\ $

{\bf Remark 1.3.} For a similar variation formula for $\eta$-invariants on odd 
dimensional manifolds with boundary, see Dai-Freed [DF].

$\ $

{\bf Remark 1.4.} For an extension of Theorem 1.2 to the  case
of families, see [DZ2].

$$\ $$

{\bf \S 2. A Riemann-Roch theorem under embedding
for Dirac operators on manifolds with boundary}

$\ $

In this section, we state a Riemann-Roch type formula for  indices of 
 Dirac type operators on 
 manifolds with boundary.
This formula will be proved in the next section and will play a key role
in our proof of the Atiyah-Patodi-Singer  index theorem in Section 4.

This section is organized as follows. In a), we describe the basic geometric
data. In b), we
state the main result  of this section, whose proof will be given in the next
section.

$\ $

{\bf a). The geometric construction of direct images under embedding
between manifolds with bounday}

$\ $

Let $Y$ be another even dimensional oriented compact
spin manifold with boundary $\partial Y$. 
Moreover, there is an
embedding $i:Y\hookrightarrow X$ such that $\partial Y\subset \partial X$,
and that $Y$ intersects transversally with $\partial X$.
 
Let $g^{T(\partial Y)}$ be the metric on $T(\partial Y)$ induced from 
$g^{T(\partial X)}$. Set
$$ {U}_\alpha' = {U_\alpha }\cap  {Y}.\eqno(2.1)$$
We can and we will assume that $\alpha$ is small enough so that
${U}_\alpha' $ is also a tubular neighborhood of $\partial Y$.
Then $ {U}_\alpha' $ carries a metric $g^{T {U}_\alpha'}$ naturally induced from
$g^{T {U}_\alpha}$.

Let
$\pi: {N}\rightarrow  {Y}$ be the normal bundle of $ {Y}$ in $ {X}$. Then 
$N_{\partial Y}= {N}|_{\partial Y}$ is the
normal bundle to $\partial Y$ in $\partial X$. 

Clearly, $\dim  {N}=\dim X-\dim Y$ is  even. Furthermore, since $TX$, $TY$ 
are oriented and spin, $ {N}$ is also oriented and spin.

Let $g^{T {X}}$ be a metric on $T {X}$ 
such that its restriction on $ {U}_\alpha$
is $g^{T {U}_\alpha }$. 
 Let $g^{T {Y}}$ be the restriction of $g^{T {X}}$ on
$ {Y}$. For simplicity, we can and we will assume that the embedding 
$ {i}:( {Y}, g^{T {Y}})\hookrightarrow ( {X},g^{T {X}})$ 
is totally geodesic. We identify $ {N}$ with
the orthogonal completement of $T {Y}$ in $(T {X})|_{ {Y}}$. 
Let $g^{ {N}}$ be the  metric
on $ {N}$ restricted from $g^{(T {X})|_{ {Y}}}$. 
Let $P^{T {Y}}$ (resp. $P^{ {N}}$) be the orthogonal
projection from $(T {X})|_{ {Y}}$ to $T {Y}$ (resp. $ {N}$) 
with respect to $g^{(T {X})|_{ {Y}}}$. Then 
$P^{T {Y}} {i}^* \nabla^{T {X}} P^{T {Y}}$,
where $\nabla^{T {X}} $ is the Levi-Civita connection of $g^{T {X}}$,
 is the Levi-Civita connection $\nabla^{T {Y}}$ of $g^{T {Y}}$ and
one has the orthogonal splitting
$$ {i}^*\nabla^{T {X}}=\nabla^{T {Y}} \oplus \nabla^{ {N}} , \eqno (2.2)$$
where $\nabla^{ {N}}=P^{ {N}} {i}^*\nabla^{T {X}}P^{ {N}}$
 is the  induced Euclidean connection
on $ {N}$.

Let $S(T {X})=S_+(T {X})\oplus S_-(T {X})$ 
(resp. $S(T {Y})=S_+(T {Y})\oplus S_-(T {Y})$,
$S( {N})=S_+( {N})\oplus S_-( {N})$) 
be the ${\bf Z}_2$-graded Hermitian vector
bundle of $(T {X},g^{T {X}})$ (resp. $(T {Y},g^{T {Y}})$, 
$( {N},g^{ {N}})$) spinors.
Then one has
  $$S\left(T {X}\right)|_{ {Y}}=S\left(T {Y}\right)
 \hat{\otimes} S\left( {N}\right).\eqno (2.3)$$
The connections $\nabla^{T {X}}$, $\nabla^{T {Y}}$, $\nabla^{ {N}}$
 lift to unitary connections
on $\nabla^{S(T {X})}$, $\nabla^{S(T {Y})}$, $\nabla^{S( {N})}$, 
$\nabla^{S^*( {N})}$ on
$S(T {X})$, $S(T {Y})$, $S( {N})$, $S^*( {N})$ 
respectively, preserving the corresponding
${\bf Z}_2$-gradings.

Let $\pi_\alpha:  {U}_\alpha  =[0,\alpha)\times \partial X\rightarrow \partial X$
(resp. $\pi_\alpha': {U}_\alpha 
=[0,\alpha)\times \partial Y\rightarrow \partial Y$)
denote the projection from $ {U}_\alpha $ (resp. 
$ {U}_\alpha' $) to the boundary of $X$ (resp. $Y$).

Let $\xi=\xi_+\oplus \xi_-$ be a ${\bf Z}_2$-graded complex vector bundle over 
$ {X}$ such that $\xi|_{ {U}_\alpha  }=\pi_\alpha^*(\xi|_{\partial X})$.
Let $g^\xi$ be a Hermitian metric on $\xi$ such that 
such that $g^\xi|_{ {U}_\alpha  }=\pi_\alpha^*(g^\xi|_{\partial X})$ and
that $\xi_+$ and $\xi_-$ are orthogonal to each other with respect to  $g^\xi$. 

Let $V\in \Gamma({\rm End}^{\rm odd}(\xi)) $ be a self-adjoint element 
such that 
$$V|_{ {U}_\alpha  }=\pi_\alpha^*\left(V|_{\partial X}\right).\eqno(2.4)$$
We assume that
$V$ is invertible on $ {X}\setminus  {Y}$, and that on $ {Y}$, 
$\ker V$ has locally constant
nonzero dimension, so that $\ker V$ is a nonzero smooth ${\bf Z}_2$-graded vector
subbundle of $\xi|_{ {Y}}$. 
Let $g^{\ker V}$ be the metric on  $\ker V$ induced by the 
metric $g^\xi|_{ {Y}}$. Let $P^{\ker V} $ be the orthogonal projection 
from $\xi|_{ {Y}}$ on $\ker V$. 

If $y\in  {Y}$, $U\in T_y {X}$, 
let $\partial _UV(y)$ be the derivative of $V$ with repsect
to $U$ in any given smooth trivialization of $\xi$ near $y\in  {X}$.
 One then verifies that
$P^{\ker V} \partial _UV(y)P^{\ker V} $ does not depend on the trivialization, and
only depends on the image $Z$ of $U\in T_y {X}$ in $ {N}_y$. 
{}From now on, we will write
$\dot{\partial}_ZV(y)$ instead of $P^{\ker V} \partial _UV(y)P^{\ker V} $. Then
one verifies easily that $\dot{\partial}_ZV(y) $ is a self-adjoint element of
${\rm End}^{\rm odd}(\ker V_y)$. 

If $Z\in  {N}$, let $\tilde{c}(Z)\in {\rm End}(S^*( {N})) $
 be the transpose of $c(Z)$ 
acting on $S( {N})$. Let $\tau^{ {N}*}\in {\rm End}(S^*( {N})) $ 
be the transpose of $\tau^{ {N}}$ defining the
 ${\bf Z}_2$-grading of $S( {N})=S_+( {N})\oplus S_-( {N})$.

Let $\mu$ be a complex vector bundle over $Y$
such that $\mu|_{ {U}_\alpha'  }=\pi_\alpha'^*(\mu|_{\partial Y})$, 
equipped with a Hermitian metric $g^\mu$
such that $g^\mu|_{ {U}_\alpha ' }=\pi_\alpha'^*(g^\mu|_{\partial Y})$.
We equip $S^*( {N})\otimes \mu$ the tensor product metric 
$g^{S^*( {N})\otimes \mu} $. 
Also, we extend an endomorphism of $S^*( {N})$ to that of 
$S^*( {N})\otimes \mu$ by acting
as identity on $\mu$. We now make the fundamental assumption that
over the total space of $ {N}$, we have the identification
$$\left( \pi^*\ker V,\pi^*g^{\ker V},\dot{\partial}_ZV(y) \right)=
\left( \pi^*\left(S^*( {N})\otimes \mu\right),
\pi^*\left(g^{S^*( {N})\otimes \mu} \right),
\tilde{c}(Z)\tau^{N*} \right).\eqno (2.5)$$

Let $\nabla^\mu$ be a Hermitian connection on $\mu$ which is of product nature
near the boundary. Let $\nabla^{S^*(N)\otimes \mu}$
be the Hermitian connection on $S^*(N)\otimes \mu$ obtained from the tensor product
of $\nabla^{S^*(N)}$ and $\nabla^\mu$.

Let $\nabla^\xi=\nabla^{\xi_+}\oplus \nabla^{\xi_-}$ be a unitary connection on 
$\xi=\xi_+\oplus \xi_-$, which preserves the ${\bf Z}_2$-grading of $\xi$ and is of
product nature near the boundary. Let
$\nabla^{\ker V}$ be the unitary connection on $\ker V$ given by
$$\nabla^{\ker V}=P^{\ker V}\nabla^\xi|_YP^{\ker V}.\eqno (2.6)$$
We then make the assumption that under the identification (2.5), we also have the
identification of connections
$$\nabla^{\ker V}=\nabla^{S^*(N)\otimes \mu}.\eqno (2.7)$$

One easily verifies that there always exists a connection $\nabla^\xi$ such that
(2.7) holds.

$\ $

{\bf Remark 2.1.} By using a well-known construction of Atiyah-Hirzebruch [AH], one
verifies easily that given metrics $g^\mu$ and $g^N$ on $\mu$ and $N$, there exist
$\xi=\xi_+\oplus \xi_-$, $g^\xi=g^{\xi_+}\oplus g^{\xi_-}$ and $V$ taken as before, 
such that (2.5) holds (Compare with [BZ, Remark 1.1]). 
In particular, $\xi_+ -\xi_-$ is a representative of
the direct image $i_!\mu\in K(X)$ of $\mu\in K(Y)$ (cf. [LM]).

$\ $

{\bf Remark 2.2.} As an easy but important observation, we note that the restriction
of the identifications (2.5), (2.7) on the boundary takes forms of exactly
the same nature (Compare with [BZ, Sect. 2b)]). In what follows, whenever such an
identification on the boundary will be considered, we will simply use a subscript
and/or superscript `$\partial$'
to indicate the restriction, when there will be no confusion
from the context.

$\ $

{\bf b). A Riemann-Roch theorem under emdedding for Dirac type 
operators on manifolds with boundary}

$\ $

We continue the discussions in a). 

Let $D^\xi=D^{\xi_+}+D^{\xi_-}$, $D^\mu$ be the Dirac operators defined as in (1.3).
We consider a Dirac type operator $D_X$ acting on $\Gamma(S(TX)\hat{\otimes} \xi)$ 
such that
$E_X=D_X-D^\xi$ is an odd endomorphism of  
$S(TX)\hat{\otimes}\xi$.

{}From (2.3), (2.5) and (A.1), one finds 
$$S(TX)|_Y \hat{\otimes}( \ker V) =S(TY) \hat{\otimes}S(N) 
\hat{\otimes}S^*(N)\otimes \mu
=(S(TY)\otimes \mu ) \hat{\otimes}\wedge^*(N^*) ,\eqno(2.8)$$
where $\wedge^*(N^*)$ is the exterior algebra bundle of $N^*$ over $Y$. 
Let $p$ be the orthogonal projection from $S(TX)|_Y {\otimes}( \ker V) $ to
$S(TY)\otimes \mu$ which maps as zero on each 
$S(TY)\otimes \mu \otimes\wedge^i(N^*),$ $i\geq 1 $.

Let $E_Y\in {\rm End}^{\rm odd}(S(TY)\otimes \mu)$ be defined by
$$E_Y=p(E_X)|_Yp.\eqno (2.9)$$

Let $D_Y$ be the Dirac type operator
$$D_Y=D^\mu+E_Y.\eqno (2.10)$$

Let $B_X$ (resp. $B_Y$) be the induced boundary operator from $D_X$ (resp. $D_Y$)
in the sense of (1.4).

The following assumption is essential for this section.

$\ $

{\bf Assumption 2.3.} {\it The operator $B_Y$ has no zero eigenvalue.}

$\ $

For any $T\in {\bf R}$, let 
$D_T:\Gamma(S(TX) \hat{{\otimes}}
{\xi}) \rightarrow \Gamma(S(TX) \hat{{\otimes}}{\xi}) $ 
be the operator defined by
$$D_T=D_X+TV,\eqno (2.11)$$
where $V\in {\rm End}(\xi)$ extends as an action on $S(TX)\hat{\otimes}\xi$
by $1\hat{\otimes}V$, etc.
Then by (1.4), its induced boundary operator $B_T$ is given by 
$$B_T=B_X-Tc\left({\partial\over\partial r}\right)V|_{\partial X}.\eqno (2.12)$$

Let $D_{T,+} $ be the restriction of $D_T$ on 
$\Gamma(S_+(TX)\otimes \xi_+\oplus S_-(TX)\otimes \xi_-)$. Let $B_{T,+}$ the associated
boundary operator and $P_{T,+,\geq 0} $ the Atiyah-Patodi-Singer projection  
associated to $B_{T,+}$.

We can now state the main result of this section as follows, 
whose proof will be given in the next section.

$\ $

{\bf Theorem 2.4.} {\it Under the Assumption 2.3, there exists $T_0>0$ such that
for any $T\geq T_0$,
$$ {\rm ind} (D_{T,+} ,P_{T,+,\geq 0} ) ={\rm ind} (D_{Y,+}, P_{Y,+,\geq 0} ).
\eqno (2.13)$$}

$\ $

{\bf \S 3. Proof of Theorem 2.4}

$\ $

The purpose of this section is to prove the Riemann-Roch property, Theorem 2.4,
for the index of boundary value problems. The proof
we described in [DZ1, Sect. 2] relies on Cheeger's cone method. Here, we will 
give a more direct proof without passing to manifolds with cone-like
singularity. We thus avoid the heat kernel analysis on cylinders and/or cones
completely.

The methods and techniques developed by Bismut and Lebeau [BL, Sects. 8, 9]
 will play an essential
role in this section. In fact, what we will do may be thought of as extensions of
the Bismut-Lebeau method to manifolds with boundary.

This section is organized as follows. In a), we construct a natural
embedding from the space of sections over $Y$ into the space of  sections over $X$. 
In b), we decompose the 
total Dirac operator on $X$ to a sum of four operators according to this 
embedding and introduce a suitable deformation
of the Dirac type operators as well as their associated boundary operators. 
In c), we prove the elliptic estimates for the deformed
operators on the boundary. In d), we prove the Fredholm property of
the Atiyah-Patodi-Singer type boundary problem for the deformed 
operators introduced in b). In e), we complete the proof of Theorem 2.4.

Throughout the rest of the paper, we will make the same assumptions 
and use the same
notation as in Section 2.

$\ $

{\bf a). An embedding mapping sections over 
$Y$ (resp. $\partial Y$) to sections over $X$ (resp. $\partial X$)}

$\ $

For any $\gamma \geq 0$, let ${\rm E}^\gamma$
(resp. ${\rm E}_{\partial X}^\gamma$, ${\rm F}^\gamma$, ${\rm F}_{\partial Y}^\gamma$) be the set of sections of
$S(TX)\hat{\otimes}\xi$ over $X$ (resp. $(S(TX)\hat{\otimes}\xi)|_{\partial X}$
over $\partial X$, $S(TY){\otimes}\mu$ over $Y$, $(S(TY){\otimes}\mu)|_{\partial Y}$ over $\partial Y$) which lie in the $\gamma^{\rm th}$ Sobolev space.

Following [BL, Sect. 8g)], for any $y\in Y$, $Z\in N_{y}$, let 
$t\in {\bf R} \mapsto x_t=\exp_y^{X}(tZ) \in X$
be the geodesic in $X$ with $x_0=y$, ${dx_t\over dt}|_{t=0}=Z$. For $ \varepsilon>0$,
set ${\cal B}_{\varepsilon } =\{ Z\in N:|Z|<\varepsilon \}$. Since $X, \ Y $ are
compact, there exists $\varepsilon_0 >0$ such that for 
$0<\varepsilon <\varepsilon_0 $, the map 
$(y,Z)\in N\mapsto \exp_y^{X}(Z) \in X$ is 
a diffeomorphism from ${\cal B}_{\varepsilon } $ onto a tubular neighborhood 
${\cal U}_{\varepsilon } $ of $Y $ in $X$. From now on, we identify 
${\cal B}_{\varepsilon } $ with ${\cal U}_{\varepsilon } $ 
and use the notation $(y,Z)$
instead of $\exp_y^{X}(Z) $. In particular, we identify $y\in Y $ with
$(y,0)\in N$.

Let $dv_{N}$ be the volume form of the fibers in $N$. Then 
$dv_{Y }(y) dv_{N}(Z) $ is a natural volume form on the total space of $N$.
Let $k(y,Z)$ be the smooth positive function on ${\cal B}_{\varepsilon_0 } $ 
defined by  
$$dv_{X}(y,Z)=k(y,Z)dv_{Y }(y) dv_{N}(Z).\eqno(3.1)  $$ 
The function $k$
has a positive lower bound on ${\cal B}_{\varepsilon_0 /2} $. Also, $k(y,0)=1$. 

Now for any $x=(y,Z)\in U_{\varepsilon_0 } $, we identify
$(S(TX)\otimes \xi)   _{x}$  
with $(S(TX)\otimes \xi)   _{y}$ 
by parallel transport with respect to $\nabla^{{S(TX)\otimes \xi}   }$ 
along the geodesic 
$t\mapsto (y,tZ)$. Clearly, this identification preserves 
the ${\bf Z}_2$-grading
of $S(TX)\hat{\otimes} \xi$.  

Take $\varepsilon \in (0,{\varepsilon _0/2}]$. Let 
$\rho :{\bf R}\rightarrow [0,1]$ be a smooth function such that
$\rho (a)=1$ if $a\leq 1/2$, while $\rho (a)=0$ if $a\geq 1$. 
For $Z\in N$, set $\rho_{\varepsilon }(Z)=\rho (|Z|/\varepsilon )$.

For  $T>0$, $y\in Y$, set 
$$\alpha_T(y)=\int_{N_{y}}\exp \left(-T|Z|^2\right) 
\rho_{\varepsilon }^2(Z) dv_{N}(Z).\eqno (3.2)$$

{\bf Definition 3.1.} {\it For any $T>0$, $\mu\geq 0$, let 
$J_T:{\rm F}^\mu\rightarrow {\rm E}^\mu$ be defined by}
$$J_T:s \mapsto 
k^{-1/2}\alpha_T^{-1/2} \rho_{\varepsilon }(Z)\exp \left(-
{T|Z|^2\over 2} \right)
s.\eqno(3.3)$$

One verifies easily that $J_T$ is well-defined. In particular,
it  induces an isometric 
embedding $J_T:{\rm F}^0\rightarrow {\rm E}^0$. 

Furthermore, one verifies that (3.3) also induces for any $T>0$, $\mu\geq 0$, 
an embedding
$$J_{T,\partial} :{\rm F}_{\partial Y}^\gamma \rightarrow
{\rm E}_{\partial X}^\gamma ,\eqno (3.4)$$
and that $J_{T,\partial}:{\rm F}^0_{\partial Y}\rightarrow {\rm E}^0_{\partial X}$
is an isometric embedding.

$\ $

{\bf b). A decomposition of  Dirac type operators under consideration
and the associated deformation}

$\ $

For any $T>0$, let
${\rm E}^0_{T}$ (resp. ${\rm E}^0_{T,\partial X}$) denote the image of
${\rm F}^0$ (resp. ${\rm F}^0_{\partial Y}$) under $J_T$ (resp. $J_{T,\partial}$).
Let ${\rm E}^{0,\perp}_{T}$ (resp. ${\rm E}^{0,\perp}_{T,\partial X}$) be the 
orthogonal completement of ${\rm E}^0_{T}$ (resp. ${\rm E}^0_{T,\partial X}$) 
in ${\rm E}^0$ (resp. ${\rm E}^0_{\partial X}$).
Let $p_T$, $p_T^\perp$ (resp. $p_{T,\partial X}$, $p_{T,\partial X}^\perp$)
be the orthogonal projections from ${\rm E}^0$ (resp. ${\rm E}^0_{\partial X}$) to
${\rm E}^0_{T}$, ${\rm E}^{0,\perp}_{T}$ (resp. ${\rm E}^0_{T,\partial X}$,
${\rm E}^{0,\perp}_{T,\partial X}$) respectively.

Recall that the Dirac operators $D_T$, $B_T$ have been defined in (2.11),
(2.12). We now decompose $D_T$, $B_T$ to
$$D_T=\sum_{i=1}^4 D_{T,i}\, , \ \ \ \ B_T=\sum_{i=1}^4 B_{T,i} \eqno(3.5)$$
respectively, where
$$D_{T,1}=p_TD_{T} p_T,\ \ D_{T,2}=p_TD_{T}
p_{T}^\perp,$$
$$D_{T,3}=p_{T}^\perp D_{T}p_{T},\ \ 
D_{T,4}=p_{T}^\perp D_{T}p_{T}^\perp.\eqno (3.6)$$
and
$$B_{T,1}=p_{T,\partial X}B_{T} p_{T,\partial X},\ \ B_{T,2}=p_{T,\partial X}
B_{T}p_{T,\partial X}^\perp,$$
$$B_{T,3}=p_{T,\partial X}^\perp B_{T}p_{T,\partial X},\ \ 
B_{T,4}=p_{T,\partial X}^\perp B_{T}p_{T,\partial X}^\perp.\eqno (3.7)$$

We now introduce a deformation of $D_T$ (resp. $B_T$) according to the decomposition
(3.6) (resp. (3.7)).

$\ $

{\bf Definition 3.2.} {\it For any $T>0$, $u\in [0,1]$, set
$$D_T(u)=D_{T,1}+ D_{T,4} +u\left( D_{T,2}+D_{T,3}\right) ,\eqno (3.8)$$
 
$$B_T(u)=B_{T,1}+ B_{T,4} +u\left( B_{T,2}+B_{T,3}\right). \eqno (3.8)'$$}

One verifies easily that $B_T(u)$ is the boundary operator associated to
$D_T(u)$ in the sense of (1.4).

$\ $

{\bf c). Elliptic estimates for $B_T(u)$}

$\ $

The purpose of this subsection is to show that the operators $B_T(u)$
verify the elliptic estimates satisfied by the usual elliptic differential
operators, when $T$ is large enough.

In fact, by
the geometric assumptions in Section 2a), when restricted to the boundary
(see in particular Remark 2.2), 
 as well as Theorem A.3, one can proceed exactly as
in [BL, Sect. 8, 9] and [BZ] to show that the following estimates for
$B_{T,i}$, $1\leq i\leq 4$, hold.

Recall that the construction of $J_T$ depends on a parameter $\varepsilon >0$.

$\ $

{\bf Proposition 3.3.} {\it There exist $\varepsilon >0$ such that
(a). as $T\rightarrow +\infty$,
$$J_{T,\partial }^{-1}B_{T,1}J_{T,\partial}
=B_Y +O\left({1\over \sqrt{T}}\right):
\Gamma \left((S(TY)\otimes \mu)|_Y\right) \longrightarrow
\Gamma \left((S(TY)\otimes \mu)|_Y\right) ;\eqno (3.9)$$
(b). there exist $C_{1}>0$, $C_{2}>0$ and $T_{0}>0$ such that for any 
$T\geq T_{0}$, 
any $s\in {\rm E}^{1,\perp }_{T,\partial X}={\rm E}^{0,\perp}_{T,\partial X} \cap
{\rm E}^{1}_{\partial X} $, $s'\in  {\rm E}^{1}_{T,\partial X} =
{\rm E}^{0}_{T,\partial X} \cap {\rm E}^{1}_{\partial X} $, then
$$\left\|B_{T,2}s\right\|_0\leq C_{1}
\left({\|s\|_1\over \sqrt{T}}+\|s\|_0\right),$$
$$\left\|B_{T,3}s'\right\|_0\leq C_{1}\left({\|s'\|_1\over \sqrt{T}}+\|s'\|_0
\right)\eqno (3.10) $$
and
$$\left\|B_{T,4}s\right\|_0\geq C_{2}
\left(\|s\|_1+ \sqrt{T}\|s\|_0\right).\eqno (3.11) $$}

{}From here, one obtains the following estimates for $B_T(u)$, 
which says that $ B_T(u)$ is a `small' perturbation of $B_T$, when $T$ is very large. 
Thus, it can be regarded as an elliptic estimate for $B_T(u)$.

$\ $

{\bf Proposition 3.4.} {\it There exist $C>0$ and $T_0>0$ such that for
any $u\in [0,1]$, $T\geq T_0$
 and $s\in {\rm E}^1_{\partial X}$, the following inequality holds,}
$$\left\| B_{T}s - B_{T}(u)s\right\|_0
\leq C\left({\left\|B_{T}s 
\right\|_0\over \sqrt{T}}+\|s\|_0\right).\eqno (3.12)$$

{\it Proof.} By the definitions of $B_T$ and $B_T(u)$, one has
$$B_Ts-B_T(u)s =(1-u)(B_{T,2}s+B_{T,3}s).\eqno (3.13)$$

{}From (3.10) and (3.13), one gets that for $u\in [0,1]$, $T\geq T_0$ 
with $T_0>0$ be as in  Proposition 3.3, one has 
$$\|B_Ts-B_T(u)s \|_0 \leq \sqrt{2}C_{1}\left({\|s\|_1\over \sqrt{T}}+\|s\|_0
\right).\eqno (3.14) $$

Now one verifies easily that the super commutator 
$[B_X,c({\partial\over \partial r})V|_{\partial X}]$ is of zeroth order.
Thus one deduces from (2.12) and the standard estimates for elliptic operators
that there exist  positive constants 
$A$, $C_3$, $C_4$ such that 
$$\| B_Ts\|_0^2\geq \| B_Xs\|_0^2- TA\| s\|_0^2$$
$$\geq C_3\| s\|_1^2-C_4\| s\|_0^2-TA\| s\|_0^2.\eqno(3.15)$$

It follows then that there exist constants $C_5>0$, $C_6>0$ such that
$$\| B_Ts\|_0\geq C_5\| s\|_1-C_6\sqrt{T}\| s\|_0.\eqno(3.16)$$

{}From  (3.14) and (3.16), one gets (3.12). Q.E.D.

$\ $

Since for any $T>0$, $B_T$ is a self-adjoint elliptic differential operator, Proposition 3.4 and the standard elliptic method enable one to deduce 
that when $T\geq \max \{T_0, 4C^2\}$, each
$B_T(u)$, for $u\in [0,1]$, is self-adjoint and 
has discrete eigenvalues with finite multiplicity.
Let $P_T(u)$ denote the Atiyah-Patodi-Singer
projection associated to $B_T(u)$. In the next subsection, we will show that
the  boundary valued problems
$(D_T(u), P_T(u))$, $u\in [0,1]$, are elliptic when $T$ is large enough.

$\ $

{\bf d). The Fredholm property of the boundary  problems
$(D_T(u), P_T(u))$}

$\ $

We continue the discussion in the previous subsection. In particular, we
assume that $T\geq \max \{T_0, 4C^2\}$ so that each $B_T(u)$, $u\in [0,1]$,  
is formally self-adjoint with discrete eigenvalues of finite multiplicity.

Set, for any $T\geq \max \{T_0, 4C^2\}$ and $u\in [0,1]$,
$${\bf E}^1_T(u) =\left\{ s\in {\rm E}^1:P_T(u)\left( s|_{\partial X}\right)
=0\right\} .\eqno(3.17)$$
Let 
$$D_{T,APS}(u) :{\bf E}^1_T(u) \longrightarrow {\rm E}^0\eqno(3.18)$$
be the uniquely determined extension of $D_T(u)$.

The main result of this subsection can be stated as follows.

$\ $

{\bf Proposition 3.5.} {\it There exists $T_1>0$ such that for any $u\in [0,1]$ and
$T\geq T_1$,
$ D_{T,APS}(u) $ is a Fredholm operator.}

{\it Proof.} By standard elliptic methods (cf. [BoW, Chap. 20]), in order to get
Proposition 3.5, it suffices to prove the following result.

$\ $

{\bf Proposition 3.6.} {\it There exist $T_1>0$, $C_7>0$, $C_8>0$ such that
for any $u\in [0,1]$,  $T\geq T_1$ and $s\in {\bf E}^1_T(u) $, one has
$$\| D_T(u)s\|_0\geq  C_7\|s\|_1 -C_8
\sqrt{T}\|s\|_0.\eqno (3.19)$$}

The rest of this subsection is devoted to a proof of Proposition 3.6.

$\ $

We decompose $X$ into two parts, the interior and the boundary region:
$$X=\left( X\setminus U_{\alpha /3}\right)\cup U_{2\alpha /3}.\eqno (3.20)$$

Our proof of Proposition 3.6 consists of three steps, corresponding 
to the interior, the boundary region and the transition region.

$\ $

{\it Step 1. The case where $s$ is supported in $ X\setminus U_{\alpha /3}$:}

Since $\alpha /3>0$, using the geometric assumptions in
Section 2a), formulas (2.9), (2.10),
Theorem A.3 and proceeding as in [BL, Sects. 8, 9] one obtains 
the following estimates.

$\ $

{\bf Lemma 3.7.} {\it There exists $\varepsilon >0$ such that
(a). as $T\rightarrow +\infty$,
$$J_{T }^{-1}D_{T,1}J_{T}
=D_Y +O\left({1\over \sqrt{T}}\right):
\Gamma \left(S(TY)\otimes \mu\right) \longrightarrow
\Gamma \left(S(TY)\otimes \mu\right) ;\eqno (3.21)$$
(b). there exist $C_{9}>0$, $C_{10}>0$ and $T_{2}>0$ such that for any 
$T\geq T_{2}$, 
any $s\in {\rm E}^{1,\perp }_{T}={\rm E}^{0,\perp }_{T} \cap {\rm E}^1 $, 
$s'\in  {\rm E}^{1}_{T} ={\rm E}^{0}_{T} \cap {\rm E}^1$
with ${\rm Supp} (|s|+|s'|)\subset X\setminus U_{\alpha /3}$, 
$$\left\|D_{T,2}s\right\|_0\leq C_{9}
\left({\|s\|_1\over \sqrt{T}}+\|s\|_0\right),$$
$$\left\|D_{T,3}s'\right\|_0\leq C_{9}\left({\|s'\|_1\over \sqrt{T}}+\|s'\|_0
\right)\eqno (3.22) $$
and
$$\left\|D_{T,4}s\right\|_0\geq C_{10}
\left(\|s\|_1+ \sqrt{T}\|s\|_0\right).\eqno (3.23) $$}

Now, from (3.21), together with the standard elliptic estimates for $D_Y$
on $Y\setminus U_{\alpha/3}'$ as well as an obvious analogue
of [BL, (9.7)], we deduce that there
exist constants $C_{11}>0$, $C_{12}>0$ such that when $T$ is large enough,
$$\left\| D_{T,1}p_Ts\right\|_0\geq C_{11}\| p_Ts\|_1 -
C_{12}\left(1+\sqrt{T}\right)\|p_Ts\|_0.\eqno (3.24)$$

Thus, using (3.8) and (3.22)-(3.24), one deduces that for $T$ large enough
and $u\in [0,1]$, 
$$\|D_T(u)s\|_0\geq \left(C_{11}-{C_{9}\over \sqrt{T}}\right)\|p_Ts\|_1
-\left(C_9+C_{12}+C_{12}\sqrt{T}\right)\|p_Ts\|_0$$
$$+\left(C_{10}-{C_{9}\over \sqrt{T}}\right)\left\|p_T^\perp s\right\|_1
+\left(C_{10}\sqrt{T}-C_9\right)\left\| p_T^\perp s\right\|_0.\eqno (3.25)$$

Estimate (3.19) follows as a consequence.

$\ $

{\it Step 2. The case where $s$ is supported in $  U_{2\alpha /3}$ :}

The key observation in this case is that since
all the geometric data are of product
nature on $U_{\alpha}$, one can use separation of variables to split the analysis into those along the
${\partial \over \partial r}$ direction and those along the cross section $\{r\}\times \partial X$'s
with $0\leq r\leq 2\alpha /3$ on which the analysis is the same as 
on $\partial X$. In particular, by (1.4), (2.11) and (2.12) one can write on
$U_{2\alpha/3}$ that
$$D_T=c\left({\partial \over \partial r}\right)\left({\partial \over \partial r}+B_T
\right). \eqno (3.26)$$

Furthermore, by the definition of the embedding $J_T$ as well as its restriction
on $\partial X$, and thus on each $\{r\}\times \partial X$, $0\leq r\leq 2\alpha /3$,
also, one deduces from (3.26) the following formula on $  U_{2\alpha /3}$,
$$D_T(u)=c\left({\partial \over \partial r}\right)
\left({\partial \over \partial r}+B_T(u) \right),\ \ \ u\in [0,1]. \eqno (3.27)$$

One also verifies easily that $B_T(u) $
anti-commutes with $c({\partial\over \partial r})$. Thus from (3.27) one gets
$$\left(D_T(u)\right)^2=-{\partial ^2\over \partial r^2}+
\left(B_T(u)\right)^2.\eqno (3.28)$$

{}From (3.27), (3.28) and  Green's formula (cf. [BoW, Chap. 3]), one deduces 
easily that
for any $s\in \Gamma (S(TX)\otimes \xi)$ which is supported in $U_{2\alpha/3}$,
$$\|D_T(u)s\|_0^2=\int_{[0,{2\alpha\over 3}]}
\left\langle B_T(u)s,B_T(u)s\right\rangle_{\{r\}\times \partial X}dr
+\left\|{\partial s\over \partial r}\right\|_0^2 - 
\langle s,B_T(u)s\rangle_{\partial X}.
\eqno (3.29)$$

Now if $s$ also verifies the boundary condition under consideration, that is,
$$P_T(u)\left(s|_{\partial X}\right)=0,\eqno(3.30)$$
then one finds
$$\langle s,B_T(u)s\rangle_{\partial X}\leq 0.\eqno (3.31)$$

On the other hand, it is clear that one can apply the analysis in Section 3c)
to each $\{r\}\times \partial X$. Thus by (3.12), (3.16) one deduces that there
exist constants $C_{13}>0$, $C_{14}>0$ such that when $T$ is large enough,
$$\|B_T(u)s\|_{\{r\}\times \partial X,0}\geq
\left(1-{C\over \sqrt{T}}\right)\|B_Ts\|_{\{r\}\times \partial X,0}
-C\|s\|_{\{r\}\times \partial X,0}$$
$$\geq C_{13}\|s\|_{\{r\}\times \partial X,1}-
C_{14}\sqrt{T}\|s\|_{\{r\}\times \partial X,0}.\eqno (3.32)$$

{}From (3.29), (3.31) and (3.32), one deduces (3.19) easily.

$\ $

{\it Step 3. The general case:}

Now by the results in Steps 1 and 2, one can apply the gluing argument 
in [BL, pp. 115-117] to complete the proof of Proposition 3.6. Q.E.D.

$\ $

The proof of Proposition 3.5 is thus also completed. Q.E.D.

$\ $

{\bf e). Proof of Theorem 2.4.}

$\ $

We assume that $T\geq T_1$ with $T_1$ determined by Proposition 3.5. 
We will first show that when $T$ is large enough, the 
family of Fredholm operators
$ D_{T,APS}(u) $, ${0\leq u\leq 1}$, constructed in Proposition 3.5 is
a continuous family. For this, one establishes the following result.

$\ $ 

{\bf Proposition 3.8.} {\it There exists $T_2>0$ such that for any 
$T\geq T_2$, $u\in [0,1]$, the operator $B_T(u)$ is invertible.}

{\it Proof.} Recall from Assumption 2.3 that $B_Y$ is invertible.
Let $c>0$ be such that 
$${\rm Spec}(B_Y)\cup [-2c,2c]=\{ 0\}.\eqno (3.33)$$
Proposition 3.8 follows from

$\ $

{\bf Lemma 3.9.} {\it There exists $T_2>0$ such that for any
$T\geq T_2$, $u\in [0,1]$ and $s\in {\rm E}^{1}_{\partial X}  $, then
$$ \left\|B_{T}(u) s\right\|_0 \geq {3c\over 2} \|s\|_0 .\eqno (3.34)$$}

{\it Proof}.  We proceed similarly as in the proof of
[TZ, Lemma 4.7]. Write $s$ as $s=s'+s''$ with $s'\in {\rm E}^{1}_{T,\partial X}$ and
$s''\in {\rm E}^{1,\perp}_{T,\partial X} $. Then one has  
$$ \left\|B_{T}(u) s\right\|_0^2=
\left\|B_{T,1}s' +  uB_{T,2}s''\right\|^2_0  +
\left\|uB_{T,3}s' +B_{T,4}s''\right\|^2_0  ,\eqno (3.35)$$
from which it follows that for any sufficiently small $\nu>0$, one has
$$\left\|B_{T}(u) s\right\|_0 \geq {7\over 8} \left\|B_{T,1}s' + u B_{T,2}s''
\right\|_0  +
\nu \left\|uB_{T,3}s' +B_{T,4}s''\right\|_0  $$
$$\geq {7\over 8} \left\|B_{T,1}s'\right\|_0 - {7\over 8} 
 \left\|B_{T,2}s''\right\|_0  +
\nu \left\|B_{T,4}s''\right\|_0- \nu \left\|B_{T,3}s'\right\|_0 .\eqno (3.36) $$

In view of (3.33), one sees easily that
$$ \left\|J_TB_YJ_T^{-1}s'\right\|_0 \geq 2c\left\|s'\right\|_0.\eqno (3.37)$$

{}From (3.37) and Proposition 3.3a), one deduces that there exists
$C_{15}>0$ such that when $T$ is sufficiently large, one has
$$ {7\over 8} \left\|B_{T,1}s'\right\|_0 \geq {3c\over 2}\|s'\|_0+{1\over 8}
\left\|J_TB_YJ_T^{-1}s'\right\|_0 -{C_{15}\over \sqrt{T}}
\left(\left\| J_TB_YJ_T^{-1}s' \right\|_0 +\|s'\|_0\right) .\eqno (3.38)$$
{}From (3.37), (3.38) one finds that when $T$ is sufficiently large,
$${7\over 8} \left\|B_{T,1}s'\right\|_0 \geq {3c\over 2}\|s'\|_0+{1\over 16}
\left\|J_TB_YJ_T^{-1}s'\right\|_0 .\eqno (3.39)$$

On the other hand, by standard elliptic estimates  as well as an obvious analogue
of [BL, (9.7)], there
exists constant $C_{16}>0$ such that
$$ \|s'\|_1 \leq C_{16} \left(\left\|J_TB_YJ_T^{-1}s'\right\|_0 +\sqrt{T}\|s'\|_0
\right) .\eqno (3.40)$$

By (3.36)-(3.40) and Proposition 3.3b), one deduces that when $T$ is 
sufficiently large, 
$$ \left\|B_T(u) s\right\|_0 \geq {3c\over 2}\|s'\|_0 +{1\over 32}
\left\|J_TB_YJ_T^{-1}s'\right\|_0 + {c\over 16}\|s'\|_0 -{7C_1 \over 8}
\left({\|s''\|_1\over \sqrt{T}}+\|s''\|_0\right)$$
$$+\nu C_2\left(\|s''\|_1+\sqrt{T}\|s''\|_0\right) 
-\nu C_1\left({C_{16}\left\|J_TB_YJ_T^{-1}s' \right\|_0\over \sqrt{T}}+
(C_{16}+1)\|s'\|_0\right) $$
$$\geq {3c\over 2} \|s'+s''\|_0 +\left({1\over 32}- 
{ \nu C_1C_{16} \over \sqrt{T}} \right) 
\left\|J_TB_YJ_T^{-1}s'\right\|_0 + \left( {c\over 16}-\nu C_1(C_{16}+1)
 \right) \|s'\|_0 $$
$$+\left(\eta C_2 - {7C_1 \over 8 \sqrt{T} }\right)\|s''\|_1
+\left(\nu C_2\sqrt{T} - {7C_1 \over 8 }-{3c\over 2} \right)\|s''\|_0.
\eqno (3.41)$$

Now  if we choose $\nu>0$ so that one also has
$${c\over 16}-\nu C_1(C_{16}+1) \geq 0,\eqno (3.42)$$
then from (3.41) one deduces easily that when $T$ is sufficiently large, 
(3.34) holds. Q.E.D.

$\ $

The proof of Proposition 3.8 is completed.
${\rm Q.E.D.}$

$\ $

{}From Proposition 3.5 and Proposition 3.8, we see that when $T$ is large
enough, we have a continuous family of Fredholm operators
$\{ D_{T,APS}(u)\}_{0\leq u\leq 1}$. Furthermore, by Proposition 3.8 and 
Green's formula, the operators 
$ D_{T,APS}(u)$, ${0\leq u\leq 1}$, are self-adjoint.

Now let $\tau_X$ (resp. $\tau_Y$) be the ${\bf Z}_2$-grading operator of
$S(TX)\hat{\otimes}\xi$ (resp. $S(TY)\otimes \mu$). One verifies directly that
$$J_T\tau_Y=\tau_XJ_T,\eqno (3.43)$$
that is, $J_T$ preserves the ${\bf Z}_2$-gradings of $S(TY)\otimes \mu$ and
$S(TX)\hat{\otimes}\xi$. 

{}From the above discussions as well as the homotopy invariance of the index of
Fredholm operators, one gets easily that
$${\rm ind} \left(D_{T,+} ,P_{T,+,\geq 0} \right) ={\rm Tr}
\left[ \tau_X|_{\ker \left(D_{T,APS}(1)\right) }\right]
={\rm Tr}\left[ \tau_X|_{\ker \left(D_{T,APS}(0)\right)}\right] .\eqno(3.44)$$

Now let $P_{T,1}$ (resp. $P_{T,4}$) be the Atiyah-Patodi-Singer projection
associated to $B_{T,1}$ (resp. $B_{T,4}$) acting on ${\rm E}^0_{T,\partial X}$
(resp. ${\rm E}^{0,\perp}_{T,\partial X}$). Then by using Proposition 3.3 and
proceed as in Section 3d), one sees easily that the boundary problems
$(D_{T,1},P_{T,1})$ and $(D_{T,4},P_{T,4} )$ are Fredholm. Furthermore, by
(3.11), (3.23) and (3.29), one deduces that when $T$ is large enough,
$$\ker \left(D_{T,4},P_{T,4} \right)=0.\eqno (3.45)$$

On the other hand, for $T$ large enough and $u\in [0,1]$, set 
$$D_Y(u)=uD_Y+(1-u)J_T^{-1}D_{T,1}J_T ,\ \ \ 
B_Y(u)=uB_Y+(1-u)J_{T,\partial}^{-1}B_{T,1}J_{T,\partial}.\eqno(3.46)$$
{}From (3.9) one can proceed as in (3.37)-(3.39) to see that when $T$ is large enough,
$B_Y(u)$ is invertible for every $u\in [0,1]$.

Let $P_{Y}(u)$ be the Atiyha-Patodi-Singer projection associated to
$B_Y(u)$. By (3.9), (3.21) and the above discussion
one sees that when $T$ is large enough,
$(D_Y(u),P_Y(u))$, $u\in [0,1]$,  form a continuous family of 
formally self-adjoint Fredholm
boundary problems. Thus by the homotopy invariance of the index of
Fredholm operators, one gets
$${\rm Tr}\left[ \tau_Y|_{\ker \left(D_Y(0),P_Y(0)\right)}\right]
={\rm Tr}\left[ \tau_Y|_{\ker \left(D_Y(1),P_Y(1)\right)}\right]
={\rm ind} (D_{Y,+}, P_{Y,+,\geq 0} ).\eqno (3.47)$$

{}From (3.43)-(3.45) and (3.47) one finds
$${\rm ind} \left(D_{T,+} ,P_{T,+,\geq 0} \right) 
={\rm Tr}\left[ \tau_X|_{\ker \left(D_{T,1},P_{T,1} \right)}\right]
+ {\rm Tr}\left[ \tau_X|_{\ker \left(D_{T,4},P_{T,4} \right) }\right]$$
$$={\rm Tr}\left[ \tau_Y|_{\ker \left(D_Y(0),P_Y(0)\right)}\right]
={\rm ind} (D_{Y,+}, P_{Y,+,\geq 0} ),\eqno(3.48)$$
which is exactly (2.13).

The proof of Theorem 2.4 is completed. Q.E.D.

$$\ $$

{\bf \S 4. The Atiyah-Patodi-Singer index theorem for Dirac operators}

$\ $

In this section, we combine Theorems 1.2 and 2.4 with the results
in [BZ] to complete our embedding proof of the Atiyah-Patodi-Singer index theorem
[APS1, (4.3)] for Dirac operators. 

This section is organized as follows. In a), we use Theorem 2.4 to refine the
main result in [BZ] so that the mod ${\bf Z}$  term in [BZ, Theorem 2.2] can now 
be made specific in our situation. In b), we 
apply the results proved in a) to the case where
$X$ is a ball to obtain the Atiyah-Patodi-Singer index theorem for Dirac operators
on $Y$.

$\ $

{\bf a). Real embeddings and $\eta$-invariants}

$\ $

Following [BZ, (1.27)], under the geometric assumptions in Section 2a), let
$\gamma^X$ be the Chern-Simons current on $X$ defined by
$$\gamma^X=\int_0^{+\infty}{\rm Tr}_s\left[V\exp \left(
-\left(\nabla^\xi +T^{1/2}V\right)^2\right)\right]
{dT\over 2T^{1/2}}.\eqno (4.1)$$

Also, if $D$ is a formally self-adjoint Dirac type operator on a closed
 odd dimensional spin manifold, we define the reduced $\eta$-invariant to be
$$\bar{\eta}(D)={\dim (\ker D)+\eta(D)\over 2}, \eqno (4.2)$$
where $\eta(D)$ is the $\eta$-invariant of $D$ in the sense of 
Atiyah-Patodi-Singer [APS1].

Let $D^{\xi_\pm}_{+,APS}$ (resp. $D^\mu_{+,APS}$) denote the Atiyah-Patodi-Singer
boundary problems associated to $D^{\xi_\pm}_+$ (resp. $D^\mu_+$). Let 
$B^{\xi_\pm}_+$ (resp. $B^\mu_+$) be the induced Dirac operators on $\partial X$
(resp. $\partial Y$) associated to $D^{\xi_\pm}_+$ (resp. $D^\mu_+$). 

Let $R^{TX}$, $R^{TY}$ denote the curvature of $\nabla^{TX}$,
$\nabla^{TY}$ respectively.

We can now state the main result of this subsection, which has
been announced in [DZ1, Theorem 2.3], as follows.

$\ $

{\bf Theorem 4.1.} {\it The following identity holds,
$${\rm ind}\left(D^{\xi_+}_{+,APS}\right) +\bar{\eta}\left(B^{\xi_+}_+\right)
-{\rm ind}\left(D^{\xi_-}_{+,APS}\right) -\bar{\eta}\left(B^{\xi_-}_+\right)$$
$$= {\rm ind}\left(D^{\mu}_{+,APS}\right) +\bar{\eta}\left(B^{\mu}_+\right)
+\left({1\over 2\pi \sqrt{-1}}\right)^{\dim X\over 2}\int_{\partial X}
{\det}^{1/2}\left({R^{TX}/2\over \sinh \left(R^{TX}/2\right)}\right)\gamma^X.
\eqno(4.3)$$}

$\ $

{\bf Remark 4.2.} A weaker mod ${\bf Z}$ version  of Theorem 4.1
has been previously proved in [BZ, Theorem 2.2]. 

$\ $

The rest of this subsection is devoted to a proof
of Theorem 4.1 
by making precise the mod ${\bf Z}$ contribution in [BZ].

For any $T\geq 0$, let $D_{T,+,APS}^\xi$ be the Dirac type operator
$$D_+^\xi+TV:\Gamma \left( \left(S(TX)\hat{\otimes} \xi\right)_+
\right) \longrightarrow
\Gamma \left( \left(S(TX)\hat{\otimes} \xi\right)_- 
\right) \eqno(4.4)$$
verifying the Atiyah-Patodi-Singer boundary condition [APS1]. Let
$B_{T,+}^\xi$ be the associated boundary operator on $\partial X$ in the sense
of (1.4).

We start with the following result which was announced in [DZ1, Prop. 2.1].

$\ $

{\bf Proposition 4.3.} {\it The quantity 
${\rm ind}(D_{T,+,APS}^\xi) +\bar{\eta}(B_{T,+}^\xi)$ does not 
depend on $T\geq 0$.}

{\it Proof.} From Theorem 1.2 and a direct counting argument in using the
definition of the reduced $\eta$-invariant, one sees easily that 
${\rm ind}(D_{T,+,APS}^\xi) +\bar{\eta}(B_{T,+}^\xi)$ depends smoothly on $T$.
Proposition 4.3 then follows from the local variation formula 
[BC3, Theorem 2.7] of Bismut and 
Cheeger.  Q.E.D.

$\ $

{\bf Remark 4.4.} To be more precise, in [BC3, Theorem 2.7], Bismut and Cheeger
considered the operators  of the form 
$B_+^\xi \tau^\xi +TV$ 
acting on $\Gamma((S_+(TX)\otimes \xi)|_Y)$, where $\tau^\xi$ is the
${\bf Z}_2$-grading operator of $\xi$. However, one sees easily that the
map $U:\Gamma((S_+(TX)\otimes \xi)|_Y) \rightarrow 
\Gamma((S_+(TX)\otimes \xi_+)|_Y \oplus (S_-(TX)\otimes \xi_-)|_Y ) $ defines by
$U:u\otimes (v_++v_-)\mapsto u\otimes v_+
-c({\partial \over \partial t})u\otimes v_-$ is an unitary and verifies that
$U( B_+^\xi \tau^\xi +TV)U^{-1}=B_{T,+}^\xi$.
This makes it clear that one can apply the results in [BC3] and [BZ] to the 
present situation.

$\ $

We can now proceed as in [BZ]. The key observation is that the geometric assumptions
in [BZ, Sect. 1b)] correspond almost exactly to the geometric assumptions on 
$\partial X$ in the current situation, with the 
minor diffference that we here use
$-c({\partial\over \partial r})\tilde{c}(Z)\tau^{N*} $ to replace 
$\sqrt{-1}\tilde{c}(Z)$ in [BZ, (1.10)]. As a result, we will use here
Theorem A.3 in the Appendix to replace [BZ, Theorem 4.5] in obtaining the 
analogues of the
analytic results of [BZ, Theorems 3.7-3.12].

We now examine the arguments in [BZ, Sect. 3e)]. In order to get the 
required result, we must find out where the mod ${\bf Z}$ terms in [BZ, Sect. 3e)]
arise, and replace them by the exact formulas.

In fact, one finds that this integer term is given by
$$ -\sharp \left\{\lambda \in {\rm Spec}\left(B_{T_0,+}^\xi \right)
: -a_0\leq \lambda < 0\right\} \eqno(4.5) $$
with $a_0>0$ such that $B^{\mu}_+ $ has no non-zero eigenvalues in
$[-2a_0,2a_0]$. This term must be added in the  right hand side of the 
analogue  of [BZ, (3.30)], 

Secondly, one uses Proposition 4.3 instead of a direct analogue of [BZ, (3.44)].

With the help of these two observations, by proceeding as in [BZ, Sect. 3e)],
one finally finds, in our situation, the following 
refinement of a direct analogue of [BZ, (3.65)] when $T_0>0$ is large enough:
$$ -\sharp \left\{\lambda \in {\rm Spec}\left(B_{T_0,+}^\xi \right)
: -a_0\leq \lambda < 0\right\} +
{\rm ind}\left(D^{\xi}_{T_0,+,APS}\right) 
- {\rm ind}\left(D^{\xi}_{0,+,APS}\right) $$
 $$ -\bar{\eta}\left(B^{\xi}_{0,+}\right)
 +\bar{\eta}\left(B^{\mu}_+ \right)
+\left({1\over 2\pi \sqrt{-1}}\right)^{\dim X\over 2}\int_{\partial X}
{\det}^{1/2}\left({R^{TX}/2\over \sinh \left(R^{TX}/2\right)}\right)\gamma^X
=0.\eqno(4.6)$$

We now prove two lemmas which together with (4.6) will give (4.3).
The first lemma follows easily from the definitions of the operators under
consideration, the definition of the reduced $\eta$-invariant as well as
the classical Agranovi\v{c}-Dynin type formula
(cf. [BoW, Chap. 21]).

$\ $

{\bf Lemma 4.5.} {\it The following identity holds,
$$ {\rm ind}\left(D^{\xi}_{0,+,APS}\right) +\bar{\eta}\left(B^{\xi}_{0,+}\right)
={\rm ind}\left(D^{\xi_+}_{+,APS}\right) +\bar{\eta}\left(B^{\xi_+}_+\right)
-{\rm ind}\left(D^{\xi_-}_{+,APS}\right) -\bar{\eta}\left(B^{\xi_-}_+\right).
\eqno (4.7)$$}

An equivalent form of the next lemma has been announced in [DZ1, Theorem 2.2].

$\ $

{\bf Lemma 4.6.} {\it The following identity holds when $T_0>0$
is sufficiently large,
$$ -\sharp \left\{\lambda \in {\rm Spec}\left(B_{T_0,+}^\xi \right)
: -a_0\leq \lambda < 0\right\} +
{\rm ind}\left( D^{\xi}_{T_0,+,APS} \right) 
= {\rm ind}\left(D^{\mu}_{+,APS} \right) .\eqno(4.8)$$}

{\it Proof.} Let $f:X\rightarrow {\bf R}$ be a smooth function 
such that $f\equiv 1$ on $U_{\alpha /3}$ and $f\equiv 0$ outside of
$U_{2\alpha /3}$.

Let $D_{T_0,-a_0,+} $ be the Dirac type operator defined by
$$D_{T_0,-a_0,+} =D^{\xi}_{T_0,+} -a_0fc\left({\partial\over\partial r}\right)
:\Gamma \left( \left(S(TX)\hat{\otimes} \xi\right)_+
\right) \longrightarrow
\Gamma \left( \left(S(TX)\hat{\otimes} \xi\right)_- 
\right). \eqno(4.9)$$
Let $D_{Y,-a_0,+} $ be the Dirac type operator defined by 
$$D_{Y,-a_0,+} =D^{\mu}_{+} -a_0fc\left({\partial\over\partial r}\right)
:\Gamma \left( S_+(TY){\otimes} \mu
\right) \longrightarrow \Gamma \left( S(TY)_-{\otimes} \mu 
\right). \eqno(4.10)$$
Let $D_{T_0,-a_0,+,APS} $, $D_{Y,-a_0,+,APS} $ be the associated operators 
verifying the Atiyah-Patodi-Singer boundary condition [APS1].

Since $-a_0$ is not an eigenvalue of $B^{\mu}_+ $, one sees that the
Assumption 2.3 is verified by the boundary operator associated to 
$D_{Y,-a_0,+} $. Thus one can apply Theorem 2.4 to get that when $T_0$ is
sufficiently large, one has
$$ {\rm ind}\left(D_{T_0,-a_0,+,APS}  \right)=
{\rm ind}\left( D_{Y,-a_0,+,APS} \right) .\eqno(4.11)$$

Now by using Theorem 1.2 and the definition of the spectral flow [APS2], 
one verifies easily that 
$$ {\rm ind}\left(D_{T_0,-a_0,+,APS}  \right)=
{\rm ind}\left( D^{\xi}_{T_0,+,APS} \right) 
-\sharp \left\{\lambda \in {\rm Spec}\left(B_{T_0,+}^\xi \right)
: -a_0\leq \lambda < 0\right\} \eqno(4.12)$$
and that
$${\rm ind}\left( D_{Y,-a_0,+,APS} \right) =
{\rm ind}\left(D^{\mu}_{+,APS} \right) .\eqno(4.13)$$

{}From (4.11)-(4.13), one gets (4.8). Q.E.D.

$\ $

{}From (4.6)-(4.8), one gets (4.3). The proof of Theorem 4.1 is now completed. Q.E.D.

$\ $

{\bf b). A proof of the Atiyah-Patodi-Singer index theorem}

$\ $

We first state an easy consequence of a direct analogue in our
situation of [BZ, Theorem 1.4].

$\ $

{\bf Lemma 4.7.} {\it The following identity holds,
$$\left({1\over 2\pi\sqrt{-1}}\right)^{\dim X\over 2}\int_X
{\det}^{1/2}\left({R^{TX}/2\over 
\sinh \left(R^{TX}/2\right)}\right)
{\rm Tr}_s\left[ \exp\left({-\left(\nabla^\xi\right)^2}\right)\right]$$
$$-\left({1\over 2\pi\sqrt{-1}}\right)^{\dim Y\over 2}\int_Y
{\det}^{1/2}\left({R^{TY}/2\over \sinh \left(R^{TY}/2\right)}\right)
{\rm Tr}\left[ \exp\left(-\left(\nabla^\mu\right)^2\right)\right]$$
$$= \left({1\over 2\pi \sqrt{-1}}\right)^{\dim X\over 2}\int_{\partial X}
{\det}^{1/2}\left({R^{TX}/2\over \sinh \left(R^{TX}/2\right)}\right)\gamma^X.
\eqno(4.14)$$ }

Finally, we are in a position to give a proof of the Atiyah-Patodi-Singer Thorem. 
Namely, we apply these results to the case where $X=D^{2n}$, the $2n$
dimensional ball with $n$ sufficiently large. 
That any compact manifold with boundary can be embedded into a 
large ball in such a fashion 
is an elementary result from differential topology.

Since $D^{2n}$ is contractable, both $\xi_\pm$ are topologically trivial
over $D^{2n}$. Thus one can deform the metric $g^{\xi_+}$ to $g^{\xi_-}$ by
$g(u)=(1-u)g^{\xi_+}+ug^{\xi_-}$, $0\leq u\leq 1$. One thus obtains easily a
smooth deformation of twisted Dirac operators moving from
$D^{\xi_+}$ to $D^{\xi_-}$. By using Theorem 1.2 as well as the standard
local variation formula of the $\eta$ invariants (cf. [APS2] and [BF, Sect. 2]), 
one gets easily the following identity
$$ {\rm ind}\left(D^{\xi_+}_{+,APS}\right) +\bar{\eta}\left(B^{\xi_+}_+\right)
-{\rm ind}\left(D^{\xi_-}_{+,APS}\right) -\bar{\eta}\left(B^{\xi_-}_+\right)$$
$$=\left({1\over 2\pi\sqrt{-1}}\right)^{\dim X\over 2}\int_X
{\det}^{1/2}\left({R^{TX}/2\over 
\sinh \left(R^{TX}/2\right)}\right)
{\rm Tr}_s\left[ \exp\left({-\left(\nabla^\xi\right)^2}\right)\right].
\eqno(4.15)$$

{}From (4.3), (4.14) and (4.15), one finds
$${\rm ind}\left(D^{\mu}_{+,APS} \right) =
\left({1\over 2\pi\sqrt{-1}}\right)^{\dim Y\over 2}\int_Y
{\det}^{1/2}\left({R^{TY}/2\over \sinh \left(R^{TY}/2\right)}\right)
{\rm Tr}\left[ \exp\left(-\left(\nabla^\mu\right)^2\right)\right]$$
$$-\bar{\eta}\left(B^{\mu}_+ \right),\eqno(4.16)$$
which is exactly the Atiyah-Patodi-singer index theorem [APS1, (4.3)] for 
$D^{\mu}_{+,APS} $.  

This completes our embedding proof of the Atiyah-Patodi-Singer index theorem 
[APS1, (4.3)] for Dirac operators on manifolds with boundary. Q.E.D.

$$\ $$

{\bf Appendix.  Dirac operators and harmonic oscillators}

$\ $

Let $E$ be a real oriented Euclidean vector space of even dimension.
Let $S(E)=S_+(E)\oplus S_-(E)$ be the ${\bf Z}_2$-graded Hermitian vector space 
of $E$-spinors.

If $e\in E$, let $e^*\in E^*$ corresponds to $e$ by the scalar product. Let
$c(e)$ denote the Clifford action of $e$ on $S(E)$. Let $\tilde{c}(e)$ denote
the corresponding Clifford action of $e$ on $S^*(E)=S_+^*(E)\oplus S_-^*(E)$.

Let $\tau$ be the ${\bf Z}_2$-grading operator of $S(E)$, that is,
$\tau|_{S_\pm(E)}=\pm {\rm id}_{S_\pm(E)}$. Let $\tau^*$ be the transpose of $\tau$.
Then $\sigma=\tau\otimes\tau^*$ is the ${\bf Z}_2$-grading operator on 
$\wedge(E^*)$.

Recall the identification of the ${\bf Z}_2$-graded vector spaces,
$$\wedge(E^*)\simeq S(E)\hat{\otimes} S^*(E).\eqno (A.1)$$

For any $e\in E$, 
let $c(e)$ (resp. $\tilde{c}(e)$) acts on $\wedge(E^*)$ as $c(e)\hat{\otimes}1 $
(resp. $1\hat{ \otimes} \tilde{c}(e)$). Then, under the identification (A.1),
$c(e)$, $\tilde{c}(e)$ acts on $\wedge(E^*)$ as
$$c(e)=e^*\wedge -i_e,$$
$$\tilde{c}(e)\tau^* =e^*\wedge +i_e\eqno (A.2)$$
respectively (Compare with [BZ, (4.5)]).

Let $e_1,\dots,e_{\dim E}$ be an oriented orthonormal base of $E$. Let
$e_1^*,\dots,e_{\dim E}^*$ be the dual base of $E^*$.

Let $\Gamma (\wedge(E^*))$ be the vector space of smooth sections of
$\wedge(E^*)$ over $E$.

$\ $

{\bf Definition A.1.} {\it Let $D^{\wedge(E^*)} $ be the operator acting on
$\Gamma (\wedge(E^*))$,
$$D^{\wedge(E^*)} =\sum_{i=1}^{\dim E} c(e_i)\nabla _{e_i},\eqno(A.3)$$
where $\nabla$ is the canonical flat connection acting on $\Gamma (\wedge(E^*))$.}

$\ $

Let $Z$ be the generic point of $E$. Then $\tilde{c}(Z)\tau^*$ acts on 
$\Gamma (\wedge(E^*))$.

$\ $

{\bf Proposition A.2.} {\it For any $T\in {\bf R}$, the following identity holds,
$$\left(D^{\wedge(E^*)} +T\tilde{c}(Z)\tau^* \right)^2=-\sum_{i=1}^{\dim E} 
\nabla_{e_i}^2 +T^2|Z|^2+T\left(\dim E -2\sum_{i=1}^{\dim E} i_{e_i}e_i^*\wedge \right).
\eqno (A.4)$$}

{\it Proof.} 
{}From (A.2),   one gets (A.4) by a direct calculation.  Q.E.D.

$\ $

Now one verifies easily that the lowest eigenvalue of 
$-\sum_{i=1}^{\dim E} i_{e_i}e_i^*\wedge$ is $-\dim E$ with the corresponding 
eigenspace being one dimensional and spanned by $1$. From this and from the standard
property of the harmonic oscillator, one gets

$\ $

{\bf Theorem A.3.} {\it The kernel of the operator 
$D^{\wedge(E^*)}+T\tilde{c}(Z)\tau^*  $ is one demensional and is spanned by 
$$\beta =\exp \left(-{T|Z|^2\over 2}\right).\eqno (A.5)$$
Furthermore, there exists $C>0$ such that $D^{\wedge(E^*)}+T\tilde{c}(Z)\tau^*  $ 
has no nonzero eigenvalue in $[-C\sqrt{T},C\sqrt{T}]$.}

$$\ $$

{\bf References}

$\ $

{\small 
[AH] M. F. Atiyah and F. Hirzebruch, Riemann-Roch theorems for differentiable 
manifolds. {\it Bull. Amer. Math. Soc.} 65 (1959), 276-281.

[APS1] M. F. Atiyah, V. K. Patodi and I. M. Singer, Spectral asymmetry and 
Riemannian geometry I. {\it Proc. Cambridge Philos. Soc.} 77 (1975), 43-69.

[APS2] M. F. Atiyah, V. K. Patodi and I. M. Singer, Spectral asymmetry and 
Riemannian geometry III. {\it Proc. Cambridge Philos. Soc.} 79 (1976),
71-99.

[AS] M. F. Atiyah and I. M. Singer, The index of elliptic operators I. 
{\it Ann. of Math.} 87 (1968), 484-530.

[B] J.-M. Bismut, Complex immersions and eta invariants. {\it Bull. Soc. Math.
France}

[BC1] J.-M. Bismut and J. Cheeger, Families index for manifolds with boundary: I, II. 
{\it J. Funct. Anal.} 89 (1990), 313-363; 90 (1990), 306-354.
 
[BC2] J.-M. Bismut and J. Cheeger, Remarks on the index theorem for families of 
Dirac operators on manifolds with boundary. in {\it Differential Geometry},
 B. Lawson and K. Tenenblat Eds., pp. 59-84, Longmam, 1992.

[BC3] J.-M. Bismut and J. Cheeger, $\eta$-invariants and their adiabatic
limits. {\it J. Amer. Math. Soc.} 2 (1989), 33-70.

[BF] J.-M. Bismut and D. S. Freed, The analysis of elliptic families II.
{\it Commun. Math. Phys.} 107 (1986), 103-163.

[BL] J.-M. Bismut and G. Lebeau, Complex immersions and Quillen metrics.
{\it Publ. Math. IHES.} 74 (1991), 1-297.

[BZ] J.-M. Bismut and W. Zhang, Real embeddings and eta invariants. 
{\it Math. Ann.} 295 (1993), 661-684.

[BeGV] N. Berline, E. Getzler and M. Vergne, {\it Heat Kernels and Dirac
Operators}. Grundl. Math. Wiss. Vol. 298, Springer-Verlag, 1992.

[BoW] B. Booss and K. P. Wojciechowski, {\it Elliptic Boundary Problems for
Dirac Operators.} Birkh\"auser, 1993.

[C1] J. Cheeger, On the spectral geometry of spaces with cone like singularities.
{\it Proc. Nat. Acad. Sci. USA.} 76 (1979), 2103-2106.

[C2] J. Cheeger, Spectral geometry of singular Riemannian spaces.
{\it J. Diff. Geom.} 18 (1983), 575-651.

[Ch] A. W. Chou, The Dirac operator on spaces with conical singularities and positive
scalar curvatures. {\it Trans. Amer. Math. Soc.}  28 (1985), 1-40.

[DF] X. Dai and D. S. Freed, $\eta$-invariants and determinant lines.
{\it J. Math. Phys.} 35 (1994), 5155-5194.

[DZ1] X. Dai and W. Zhang, The Atiyah-Patodi-Singer index theorem for manifolds
with boundary: a proof using embeddings. {\it C. R. Acad. Sci. Paris,
S\'erie I}, 319 (1994), 1293-1297.

[DZ2] X. Dai and W. Zhang, Higher spectral flow. {\it J. Funct. Anal.} 
157 (1998),  432-469.

[DZ3] X. Dai and W. Zhang, Families of real embeddings, $\hat{\eta}$-forms
and  index bundles for Dirac operators on manifolds with boundary. {\it To appear.}

[LM] H. B. Lawson and M.-L. Michelsohn, {\it Spin Geometry.}
 Princeton Univ. Press, 1989.

[M] R. B. Melrose, {\it The Atiyah-Patodi-Singer index theorem.} A. K. Peters Ed.,
1994.

[MP] R. B. Melrose and P. Piazza, Families of Dirac operators, boundaries and
the $b$-calculus. {\it J. Diff. Geom.} 45 (1997), 99-180.

[TZ] Y. Tian and W. Zhang, Quantization formula for symplectic manifolds with
boundary. {\it Geom. Funct. Anal.} To appear.

[Z] W. Zhang, A proof of the mod 2 index theorem of Atiyah and Singer. {\it C. R.
Acad. Sci. Paris, S\'erie I}, 316 (1993), 277-280.}

$\ $

X. D: Department of Mathematics, University of  California,
Santa Barbara, California 93106, USA

{\it E-mail: dai@math.ucsb.edu}

$\ $

W. Z: Nankai Institute of Mathematics, 
Nankai University, Tianjin 300071, People's  Republic
of  China

{\it E-mail: weiping@sun.nankai.edu.cn}

\end{document}